\documentclass[12pt]{article}

\usepackage{latexsym,amssymb}
\newcommand{\hq}{H_1(G;\mathbb Q)}
\newcommand{\hkz}{H_1(K;\mathbb Z)}
\newcommand{\hkq}{H_1(K;\mathbb Q)}
\newcommand{\z}{\mathbb Z}
\newcommand{\q}{\mathbb Q}
\newtheorem{lem}{Lemma}[section]
\newtheorem{defn}[lem]{Definition}
\newtheorem{ex}[lem]{Example}
\newtheorem{co}[lem]{Corollary}
\newtheorem{thm}[lem]{Theorem}
\newtheorem{prop}[lem]{Proposition}

\newenvironment{proof}{\textbf{Proof.}}{\newline\hspace*{\fill}{$\Box$}}

\begin{document}
\title{Finite covers of the infinite cyclic cover of a knot}
\author{J.\,O.\,Button\\
Selwyn College\\
University of Cambridge\\
Cambridge CB3 9DQ\\
U.K.\\
\texttt{jb128@dpmms.cam.ac.uk}}
\date{}
\maketitle
\begin{abstract}
We show that the commutator subgroup $G'$ of a classical knot group $G$
need not have subgroups of every finite index, but it will if $G'$ has a 
surjective homomorphism to the integers and we give an exact criterion for
that to happen. We also give an example of a smoothly knotted $S^n$ in
$S^{n+2}$ for all $n\geq 2$ whose infinite cyclic cover is not simply
connected but has no proper finite covers.
\end{abstract}
\section{Introduction}

Given a classical knot, that is a piecewise linear embedding of $S^1$ into
$S^3$, the fundamental group $G$ of its complement (which we will call a
classical knot group) contains a lot of information about the knot itself
and is a much studied object. We have that the quotient of $G$ by its
commutator subgroup $G'$ is $\z$ and this immediately implies that $G$ has
subgroups of every finite index. Using the correspondence between
subgroups and covers, we can also regard this as telling us that 
the knot complement has coverings of every finite degree. Moreover there is
the cover given by $G'$, known as the infinite cyclic covering, and this cover
will be simply connected if and only if the knot is trivial. We can ask what 
coverings of finite degree will be possessed by the infinite cyclic covering 
of a non-trivial knot. Here we split into two radically different cases
because if the knot is fibred, which happens if and only if $G'$ is
finitely generated, then $G'$ is a free group with rank equal to the
genus of the knot and so has subgroups of every finite index. However the
focus throughout this paper is when $G'$ is not finitely generated. In
this case it is known by \cite{sw1} that $G'$ can have infinitely many,
indeed uncountably many, subgroups of a particular finite index and in 
\cite{swi} a related result states that if $G'$ has infinitely many
representations into $S_3$ then it has uncountably many subgroups of every
index $n\geq 3$ (and some of index 2, although the commutator subgroup
of a knot group can only ever have finitely many subgroups of index 2).

In this paper we first give in Section 2 a condition on the Alexander
polynomial of our knot which instantly tells us exactly when the commutator
subgroup $G'$ has subgroups of index 2. As there are plenty of knots which
fail this condition, we have an answer to the question of Silver and
Williams in \cite{swi} which asks whether the commutator subgroup of a
non-fibred knot maps onto $S_n$ for all $n$ and hence has subgroups of
every finite index. However we also give an infinite family of two bridge
knots which have infinitely many representations into $S_3$ and hence
subgroups of every finite index by the result above.

Another way to conclude that a group has subgroups of every finite index
is to find a homomorphism onto $\z$. In Section 3 we show that the
commutator subgroup of a knot group possesses such a homomorphism if and
only if the Alexander polynomial of the knot has a non-trivial factor which
is monic at both ends. As Alexander polynomials are readily computed from
a diagram of the knot or a presentation of the knot group, we can work out
whether this is satisfied for any given knot.

In Section 4 we move to higher dimensional knots, which are locally flat
embeddings of $S^n$ in $S^{n+2}$ for $n\geq 2$. There we still have the
knot group $G$ and commutator subgroup $G'$ with $G/G'=\z$ so that $G'$
corresponds to the infinite cyclic covering of the knot, but here the
position regarding finite index subgroups of $G'$ can be quite different.
We show that there exist knots in each dimension such that the infinite
cyclic covering is not simply connected but has no proper coverings of 
finite degree. We do this by utilising a result of Kervaire which gives
us sufficient conditions for a finitely presented group to be a higher
dimensional knot group. We can also say that these knots are smooth, which
we conclude from Kervaire's result when $n\geq 3$ and by a theorem of
Levine when $n=2$.

\section{Examples}

Given a finitely generated group $G$, it is a standard definition that
$\beta_1(G)$ is the number of $\z$ summands in the abelianisation
$G/G'=H_1(G;\z)$. But if $G$ is not finitely generated then two
alternatives present themselves; as $H_1(G;\q)=H_1(G;\z)\otimes_{\z}\q$
is a vector space over $\q$, we let $\beta_1(G;\q)$ be the dimension of
this vector space. However a natural definition of $\beta_1(G;\z)$ would
be the maximum integer $k$ such that there is a surjective homomorphism
from $G$ to $\z^k$. As this implies that $\z^k\leq G/G'$, we have 
$\hq$ contains a linearly independent subspace of dimension $k$ and thus
$\beta_1(G;\z)\leq\beta_1(G;\q)$.

Here we are interested in groups $K=\mbox{ker }\chi$ where $\chi$ is a
surjective homomorphism from a finitely presented group $G$ to $\z$.
These groups may or may not be finitely generated and it is the infinitely
generated case in which we are most interested. However $K$ will be better
behaved than an arbitrary infinitely generated group as the 
Reidemeister-Schreier rewriting process gives us a presentation of $K$ 
and thus,
by abelianising, of $K/K'=\hkz$. This means that $K/K'$ can be regarded as a
finitely presented module over the unique factorisation domain 
$\z[t^{\pm 1}]$, where $t$ acts by conjugation on $K$ using an element of
$\chi^{-1}(1)$. In this situation
we can gain a lot of information by using the Alexander
polynomial.
\begin{defn}
Given a finitely presented module $M$ over a unique factorisation domain $R$
with $P$ an $m$ by $n$ presentation matrix for $M$, we define the Alexander
polynomial $\Delta_M\in R$ to be the highest common factor of the $m\times m$
minors of $P$ (where we assume that $m\leq n$ by adding zero columns if
necessary).
\end{defn}
This is independent of the presentation matrix and is defined up to
multiplication by units in $R$.

Where $K$ is above we write the Alexander polynomial as $\Delta_{G,\chi}(t)\in
\z[t^{\pm 1}]$ and it is defined up to multiplication by $t^{\pm k}$.
We can define the degree of $\Delta_{G,\chi}(t)$ as the degree
of the highest power of $t$ minus the lowest (with 0 having infinite degree).
In this case we can also consider the module $\hkq$ over $\q[t^{\pm 1}]$,
whereupon we can use the structure theorem for finitely generated modules
over a principal ideal domain to conclude that $\hkq$ is a direct sum of
cyclic modules. In particular it has a square presentation matrix whose
determinant is $\Delta_{G,\chi}(t)\in\q[t^{\pm 1}]$ (although now
$\Delta_{G,\chi}(t)$ is only defined up to multiplication by $qt^{\pm 1}$
for $q\in\q-\{0\}$)
and so we conclude 
that the degree of $\Delta_{G,\chi}$ is $\beta_1(K;\q)$. However when $K$ is
infinitely generated this may be strictly larger than $\beta_1(K;\z)$.

\begin{ex} 
\end{ex}
Let $G=\langle a,t|tat^{-1}=a^2\rangle$ with $\chi(t)=1,\chi(a)=0$ and
$K=\mbox{ker }\chi$ so that $\Delta_{G,\chi}
(t)=t-2$ with $\beta_1(K;\q)=1$. However $K$ is well known to be isomorphic
to the dyadic rationals and has no surjective homomorphism to $\z$ (for
instance $K$ has no subgroups of index 2). Hence $\beta_1(K;\z)=0$.

In the above we interpret the degree of $\Delta_{G,\chi}$ as being infinity
if and only if $\Delta_{G,\chi}=0$. We see this happens if and only if
$H_1(K;\q)$ has a non-trivial free $\q[t^{\pm 1}]$-submodule which is
equivalent to the dimension of $H_1(K;\q)$ being infinite. In this
case we can conclude by a result \cite{how} of Howie
that the finitely presented group $G$ is large, that is
it has a finite index subgroup possessing a surjective homomorphism to a
non-abelian free group, which has a range of consequences.

Given the abelianised relations obtained for $K$ which yields defining
equations for $K/K'$, in addition to considering these equations over $\q$
we can also look at them over the finite field $\z/p\z$ for $p$ a prime.
Again we find that $H_1(K;\z/p\z)$ is a finitely presented module over
$\z/p\z[t^{\pm 1}]$ and this Alexander polynomial is obtained by reducing
the standard Alexander polynomial modulo $p$. Again Howie's result gives
us largeness if this vanishes.
\begin{prop}
Given $K=\mbox{ker }\chi$ for $\chi$ a surjective homomorphism from the
finitely presented group $G$ to $\z$ then the number $r_p$ of representations
of $K$ into $\z/p\z$ is $p^{d(p)}$ and the number $n_p$ of normal subgroups
of $K$ of prime index $p$ is $(p^{d(p)}-1)/(p-1)$ where $d(p)$ is the degree
of the mod $p$ reduced Alexander polynomial $\Delta_{G,\chi}$.
\end{prop}
\begin{proof}
A normal subgroup of index $p$ in $K$ gives rise to a surjective homomorphism
from $K$ to $\z/p\z$. Now Hom$(K;\z/p\z)\cong H^1(K;\z/p\z)\cong (\z/p\z)
^{d(p)}$ and any non-zero homomorphism is surjective. Moreover if two
homomorphisms have the same kernel then they are related by an automorphism
of the quotient.
\end{proof}

This gives an immediate corollary.
\begin{co}
The group $K$ as above has infinitely many normal subgroups of prime index
$p$ if and only if the Alexander polynomial $\Delta_{G,\chi}$ is 0 modulo
$p$. Also $K$ has no normal subgroups of prime index $p$ if and only if
$\Delta_{G,\chi}$ is constant but non-zero modulo $p$. 
\end{co}

Now let $G$ be a (non-trivial)
classical knot group with $K$ equal to the commutator
subgroup. It is noted in $\cite{swi}$ that if the knot is fibred then
$K$ (which will be free of finite rank) contains subgroups of every finite 
index and it is then asked whether this is always true in the non-fibred case.
From the above we gain an immediate answer of no.
\begin{co}
The commutator subgroup of a classical knot group has no subgroups of order 2
if and only if the Alexander polynomial
\[\Delta(t)=c_n(t^n+t^{-n})+\ldots +c_1(t+t^{-1})+c_0\]
of the knot has all coefficients even except $c_0$.
\end{co}

Note that $c_0$ is always odd because $\Delta(1)=\pm 1$. Moreover, as any
symmetric polynomial $f(t)$ of even degree with $f(1)=\pm 1$ is the
Alexander polynomial of some knot (see \cite{lev65}), we have infinitely
many knots of arbitrarily high genus with this property as well as 
all knots with Alexander polynomial 1. 
 
In \cite{swi} Corollary 3.9 it is shown that
a sufficient condition for the commutator subgroup $K$ of a classical knot
group to have subgroups of every finite index is
that $K$ has infinitely many representations into the symmetric group $S_3$:
in fact it is actually shown that 
this ensures uncountably many subgroups for every index at
least 3, whereas there can only be finitely many subgroups of index two
because for a knot we have $\Delta(1)\equiv 1$ mod 2. It seems appropriate
here to show that this case actually occurs for an infinite family of
examples. 
\begin{thm}
For any positive $n$ which is 3 mod 6, the commutator subgroup $K$ of the 
2-bridge knot $(4n-1/4n-3)$ has infinitely many representations into $S_3$
and hence $K$ has subgroups of every finite index.
\end{thm}
\begin{proof}
We adopt the notation of \cite{mr} in obtaining a presentation of the
fundamental group $G$ of the complement of the 2-bridge knot $(p/q)$
where $p$ and $q$ are coprime with $p$ odd and $0<q<p$. We form the 
sequence $(k_i)$ for $1\leq i\leq p-1$ by
setting $iq=k_ip+r_i$, where $0<r_i<p$ so that $k_i$ is the quotient of
$iq$ by $p$ under Euclidean division.
We let $e_i=(-1)^{k_i}$ thus the
number of successive $e_i$s with the same sign is just the cutting sequence
of $p/q$ which in turn is closely related to its continued fraction.
Then $G$ has the presentation
\[\langle
u,v|w=v^{e_1}u^{e_2}\ldots u^{e_{p-2}}v^{e_{p-1}},uw=wv\rangle\]
for $u$ a meridian of the knot.

For our case where $(p/q)=(4n-1/4n-3)$, the $e_i$s alternate in sign
successively from $e_1=+1$ to $e_{2n-1}=+1$ and alternate again from
$e_{2n}=+1$ to $e_{4n-2}=+1$. Now on replacing $v$ with $ua$ in the relation,
we see that for the
2-bridge knot $(4n-1/4n-3)$ we obtain the presentation
\[\langle u,a|ua^nua^{-n}u^{-1}a^{n-1}u^{-1}a^{-n}\rangle\]
meaning that by the Reidemeister-Schreier rewriting process we have
\[K=\langle a_i|a_{i+1}^na_{i+2}^{-n}a_{i+1}^{n-1}
a_i^{-n}\mbox{ for }i\in\z\rangle\]
where $a_i$ is the element $u^iau^{-i}$ in $G$.
Thus our representations of $K$ into $S_3$ can be thought of as biinfinite
sequences $(a_i:i\in\z)$ of elements of $S_3$ such that the above relation
holds. Now we use the fact that $n$ is odd but divisible by 3. If both
$a_i$ and $a_{i+1}$ are the 2-cycle $(xy)$ then we require $a_{i+2}$ such
that $(xy)a_{i+2}^n(xy)=\mbox{id}$ so that $a_{i+2}$ can be either of
the 3-cycles or the identity. On taking $a_{i+2}$ to be the latter we find
on replacing $i$ with $i+1$ and then $i+2$
that we must have $a_{i+3}=a_{i+4}=(xy)$ and so we return to our initial
conditions. However if we choose $a_{i+2}=(123)$ say, we find that
$a_{i+3}^{-n}=(xy)(123)$ so $a_{i+3}$ moves on to another 2-cycle. But then
we have $a_{i+3}a_{i+4}^{-n}=\mbox{id}$ so that $a_{i+4}=a_{i+3}$.
From above we have a choice of directions and by taking
$a_{i+5}=(123)$ again we can ensure that $a_{i+6}=a_{i+7}$ is the third
2-cycle, hence on taking $a_{i+8}=(123)$ we get back to $a_{i+9}=a_{i+10}=
(xy)$ and so we have a sequence of period 9. But at any point where we
have the same 2-cycles for $a_j$ and $a_{j+1}$ we can loop through the
identity rather than sticking to the periodic orbit, thus we have
uncountably many different representations of $K$ into $S_3$.
\end{proof}

\section{Surjections to the integers}

We now vary our approach in trying to find knots whose commutator subgroup
$K$ has subgroups of every finite index. An immediate condition that implies
this is that $K$ has a surjection to $\z$. However it is not enough just to
look at the degree of the Alexander polynomial as Example 2.2 shows.
However let us first assume that our knot group $G$ has Alexander polynomial
$\Delta(t)=a_dt^d+\ldots +a_0$ with $a_d$ and $a_0$ not equal to 0
and that the abelianisation $H_1(K;\z)$ of
$K$ is a cyclic module over $\z[t^{\pm 1}]$ which always happens if $G$ has
a 2-generator 1-relator presentation. We then have that $H_1(K;\z)=\z
[t^{\pm 1}]/(\Delta)$ as a $\z[t^{\pm 1}]$-module and therefore as an
abelian group $A$ it has generators $x_i$, where $x_i$ corresponds to $t^i$,
and presentation
\[A=\langle x_i|a_nx_{n+i}+\ldots +a_1x_{i+1}+a_0x_i=0
\mbox{ for }i\in\z\rangle.\]
Therefore the way here to think about homomorphisms from $A$ to $\z$ is that 
they are biinfinite recurrence sequences ($x_i:i\in\z$) satisfying the above
homogenous linear recurrence relation. Therefore if $\Delta$ is non-constant
and both $a_d$ and $a_0$ are $\pm 1$ then we obtain a non-trivial
homomorphism (and hence a surjective one) from $A$ to $\z$ by picking
arbitrary integers not all zero for the initial conditions $x_0,\ldots ,
x_{d-1}$ and then we use the equations to find $x_d,x_{d+1},\ldots$ and
$x_{-1},x_{-2},\ldots$ which will all be integers.
Of course if the knot is fibred then the Alexander
polynomial is monic and $A$ is just $\z^d$ for $d$ the degree of $\Delta(t)$.
However it is well known that there are knots with monic and non-constant 
Alexander polynomial which are not fibred. 
Moreover it is also clear in the case of cyclic modules
that if $\Delta$ has a non-constant factor $f$ 
which is monic at both ends then we obtain a non-trivial
solution in integers to the linear recurrence relation defined by $f$, and
hence for that defined by $\Delta$, so again finding a surjective
homomorphism from $A$ to $\z$ and hence from $K$ to $\z$.

This also works if $H_1(K;\z)$ is isomorphic to a direct sum of cyclic
modules
\[\z[t^{\pm 1}]/(f_1)\oplus\ldots\oplus\z[t^{\pm 1}]/(f_m)\]
as then $\Delta=f_1\ldots f_m$. Thus if $\Delta$ has a factor $f$ which is
monic at both ends then we can assume it is irreducible and hence $f$
divides some $f_i$. We can then use the homomorphism from $H_1(K;\z)$ to
$\z[t^{\pm 1}]/(f_i)$ and then argue as above. In order to argue generally we 
follow Lemma 3.6 in \cite{swtag} which states that if $M$ is a finitely
generated $\z[t^{\pm 1}]$-module which as an abelian group is torsion free
then $M$ is module isomorphic to a finite index submodule of
\[M'=\z[t^{\pm 1}]/(\pi_1)\oplus\ldots\oplus\z[t^{\pm 1}]/(\pi_n)\]
where the polynomials $\pi_i$ are primitive with $\pi_i|\pi_{i+1}$.
\begin{thm}
If a knot in $S^3$ has a non-constant factor $f(t)$ of its Alexander polynomial
$\Delta(t)$ which is monic in $t$ and in $1/t$ then the commutator subgroup
$K$ of the knot group has a homomorphism onto $\z$ and hence $K$ has subgroups
of every finite index.
\end{thm}
\begin{proof}
On setting $M$ to be the module $H_1(K;\z)$ as before, we do have that $M$
is $\z$-torsion free because any knot group 
has a deficiency 1 presentation and hence a square 
presentation matrix $P$ with det $P=\Delta(t)$. Thus suppose we had a
non-zero $y\in M$ with $py=0$ for some prime $p$. This corresponds to a
column vector $\mu$ which is not in the image of $P$ but with a vector
$\lambda$ such that $P\lambda=p\mu$. We can pull out multiples of $p$ so
that on looking modulo $p$ we have $\lambda\not\equiv 0$ with $P\lambda
\equiv 0$. Thus det $P\equiv 0$ but Alexander polynomials of knots
satisfy $\Delta(1)=\pm 1$.

Let us identify $M$ with its image in $M'$
and take a non-zero\\ 
$m=(m_1,\ldots ,m_n)\in M$ such that $m$ is
annihilated by our factor $f(t)$ which we can take to be irreducible.
This $m$ exists because we have det $P\equiv 0$
mod $f$ so there exists a column vector $u$ with
$Pu\equiv 0\mbox{ mod }f$ but $u\not\equiv 0\mbox{ mod }f$. 
Let $m$ be defined by $Pu=fm$ so that $fm$ is in the image of $P$, which
is 0 in $M$. Hence we can take this $m$ provided it is not itself in the
image of $P$. But then if $m=Pn$ we can invert $P$ in the field of fractions
of $\z[t^{\pm 1}]$ because det $P$ is not 0 in $\z[t^{\pm 1}]$. So if
$m=Pn$  then this gives 
$u=fn$ so $u$ was 0 mod $f$ which is a contradiction.

Now take $k$ such that $m_k$ is non-zero modulo $\pi_k$ and then we must have
$f$ dividing $\pi_k$. This means that $\z[t^{\pm 1}]/(\pi_k)$ maps onto $\z$
thus $M'$ does too and then we can restrict this homomorphism to $M$, in
which case the image is also a copy of $\z$ because $M$ has finite index
in $M'$.
\end{proof}

We also have a converse to Theorem 3.1.
\begin{thm}
A knot in $S^3$ has the property that the commutator subgroup $K$ of the
knot group $G$ possesses a surjective homomorphism to $\z$ if and only if
there is a non-constant factor of its Alexander polynomial which is 
monic at both ends.
\end{thm}
\begin{proof}
We suppose that $K$ has a surjection to $\z$ and
as before we reduce to the case where the abelianisation $H_1(K;\z)$
of $K$ is a cyclic $\z[t^{\pm 1}]$-module. Setting $M=H_1(K;\z)$,we use the
fact that $M$ has finite index in the module $M'$ which is a direct sum of
cyclic modules. For any subgroup $N$ of $M$ we have that
$M'/N$ has a normal subgroup $M/N$ of finite index. So if $M/N=\z$ then
$M'/N$ is a
finitely generated infinite abelian group and hence has a surjection to $\z$,
thus so does $M'$. Next we restrict this homomorphism to each of the cyclic
modules so that there must be a $k$ with $\z[t^{\pm 1}]/(\pi_k)$ surjecting to
$\z$. Now suppose that $f$ is a non-constant irreducible factor of $\pi_k$.
We can find $m\in M$ with $m\neq 0$ but
$fm=0$ so, by the same argument as in Theorem 3.1,
we conclude that $f$ divides the Alexander
polynomial $\Delta(t)$. Hence by the comment at the beginning of this section,
we are done if we prove that whenever we have a biinfinite recurrence 
sequence $(x_n)$ lying entirely in $\z$ that satisfies the recurrence
relation defined by
\[a_dx_{n+d}+\ldots +a_0x_n=0\mbox{ for all }n\in\z\]
where $d\geq 1$ with $a_0,\ldots ,a_d\in\z$ and $a_0,a_d\neq 0$ then the
auxiliary polynomial $a_dt^d+\ldots +a_0$ has a factor that is monic at
both ends. This is established in Theorem 3.3.
\end{proof}

We did not find a proof of the following result in the literature, probably
because the standard definition of a recurrence sequence uses
a monic polynomial. The author would like to thank R.\,G.\,E.\,Pinch for
suggesting the following argument.
\begin{thm}
Given a non-trivial biinfinite recurrence sequence $(x_n:n\in\z)$ with all 
terms in $\z$ that satisfies the relation
\[a_dx_{n+d}+\ldots +a_0x_n=0\mbox{ for all }n\in\z\]
where $d\geq 1$ with $a_0,\ldots ,a_d\in\z$ and $a_0,a_d\neq 0$ then the
auxiliary polynomial $a_dt^d+\ldots +a_0$ has a non-constant factor that is 
monic at both ends.
\end{thm}
\begin{proof}
The set of all polynomials such that $(x_n)$ satisfies the corresponding
recurrence relation forms an ideal 
in $\q[t]$
which is a principal ideal domain, so we take a generator $f(t)$. On
clearing out denominators and removing any integer dividing all the
coefficients, we can assume that $f(t)$ lies in $\z[t]$ and is defined
uniquely up to sign (which we think of as the ``correct'' auxiliary
polynomial). We will treat $(x_n)$ as an ordinary recurrence relation with
$n\in\mathbb N$ and will show that all the factors of $f(t)$ are monic.
This will then imply our result in the biinfinite case by sending $n$ to
$-n$ and $f(t)$ to $t^df(t^{-1})$.

We use the theory of local fields; see \cite{cas} Chapter 10 or \cite{ser}
Chapter 1 for details. Let the roots of the auxiliary polynomial $f$ be
$\alpha_1,\ldots ,\alpha_d$. We first assume that $f$ is irreducible.
Consequently we have the solution
\[x_k=\beta_1\alpha_1^k+\ldots +\beta_d\alpha_d^k.\]
Let $K$ be the number field obtained from $\q$ by adjoining all the roots
of $f$, which will also contain all $\beta_i$ as can be seen by looking
at the Vandermonde matrix for the first $d$ terms of the sequence.
Moreover this shows that no $\beta_i$ is zero (unless $x_k$ is identically
zero) because a field automorphism $\sigma$ of $K$ which sends $\alpha_i$
to $\alpha_{\sigma(i)}$ also sends $\beta_i$ to $\beta_{\sigma(i)}$.

Assuming that $f$ is not monic means that there is a discrete
valuation $v$ on $K$ and an $\alpha_i$ with $v(\alpha_i)<0$. Let $\pi$ be
a uniformising element for $K$, so that any $k\in K-\{0\}$ can be expressed as
$k=\pi^nu$ for $n\in\z$ and $v(u)=0$, giving $v(k)=n$. 
If $A$ is the valuation ring (namely
the set of those $k\in K$ with $v(k)\geq 0$) then $u$ is a unit of $A$ if
and only if $v(u)=0$. Thus let us write $\alpha_i=\gamma_i/\pi^{e_i}$ and
$\beta_i=\delta_i/\pi^{f_i}$ where the $\gamma_i$ and the $\delta_i$ are all
units. We reorder so that $e_1\geq e_2\geq\ldots \geq e_l$ with $r$ such that
$e_1=\ldots =e_r>e_{r+1}$. Note that $e_1>0$. Also we reorder the $f_i$ so
that $f_1\geq \ldots \geq f_r$.

Now consider the element
\[\delta=\delta_1+\pi^{f_1-f_2}\delta_2+\ldots +\pi^{f_1-f_r}\delta_r\]
and take $c\in\mathbb N$ such that 
$c>v(\delta)\geq 0$. As the quotient of $A$ by
$\pi^cA$ is a finite ring, say of order $N$, we have for 
any unit $u$ of $A$ that
$u^N\equiv 1$ mod $\pi^cA$. Thus for all $k\in\mathbb N$ we obtain
\[x_{Nk}=\beta_1\alpha_1^{Nk}+\ldots +\beta_d\alpha_d^{Nk}=
\frac{\delta_1(1+\pi^ca_1)}{\pi^{e_1Nk+f_1}}+\ldots +
\frac{\delta_d(1+\pi^ca_d)}{\pi^{e_dNk+f_d}}\]
where $a_1,\ldots ,a_d\in A$.
As $c>0$ we have for $j\geq r+1$ that the $j$th term in the sum has 
valuation equal to $-(e_jNk+f_j)$. Now pick $k$ large enough to ensure that
\[e_1Nk+f_1-c>\mbox{max }\{e_jNk+f_j:r+1\leq j \leq d\}.\]
Then on putting the first $r$ terms over a common denominator to form 
$t_{Nk}$ we have
\[t_{Nk}=(\delta_1(1+\pi^ca_1)+\pi^{f_1-f_2}\delta_2(1+\pi^ca_2)
+\ldots + \pi^{f_1-f_r}\delta_r(1+\pi^ca_r))/(\pi^{e_1Nk+f_1})\]
and the numerator has the same valuation as $\delta$ so that overall we have
\[v(t_{Nk})=v(\delta)-(e_1Nk+f_1)<\mbox{min }\{-(e_jNk+f_j):
r+1\leq j\leq d\}.\]
Hence as 
\[x_{Nk}=t_{Nk}+\beta_{r+1}\alpha_{r+1}^{Nk}+\ldots +\beta_d\alpha_d^{Nk}\]
we have $v(x_{Nk})=v(t_{Nk})$. 
Consequently as soon as $k$ is large enough to make
$v(t_{Nk})<0$, we have $v(x_{Nk})<0$ so we do not have a sequence of integers.

As for the case when $f$ is reducible or even has
repeated factors, the quickest way to
deal with this case is to consider the shift map
$S((x_n))=(x_{n+1})$ which, as a linear operator on the space of rational
sequences satisfying our recurrence relation, has minimal polynomial 
equal to our auxiliary equation $f$. Therefore if $f=gh$ we can form a
new sequence $y_n$ by applying $g(S)$ to $x_n$ and this will satisfy
the recurrence defined by the polynomial $h$. Moreover as $g$ has integer
coefficients we see that $y_n$ is also an integer valued sequence. We can
then apply the above to $y_n$, thus concluding that the 
irreducible factors of $f$ are monic.
\end{proof}

\section{Higher dimensional examples}

In the previous two sections we have seen examples both of knots whose
commutator subgroups have subgroups of every finite index and those which
do not. However we can certainly say that as classical knot groups are
residually finite then so too are their commutator subgroups $K$, and any
infinite residually finite group has infinitely many finite index subgroups
of arbitrarily high index. Here we go in the other direction and look at
finitely presented groups $G$ with a homomorphism onto $\z$ where the
kernel has no proper finite index subgroups. Of course one way of
achieving this is to take a finitely presented group $H$ with this property
and then set $G$ to be the direct product $H\times\z$ (or even a semi-direct
product). However we wish to come up with examples which look more like
classical knot groups; in particular we require that our group $G$ has
$H_1(G;\z)=\z$, $H_2(G;\z)=0$, and weight 1 which means that there exists
an element whose normal closure is all of $G$ (such as a meridian in the case
of a classical knot group). Let us for this purpose
use presentations of deficiency 1, and in particular
2-generator 1-relator presentations $\langle x,y|r\rangle$. We can
always make a change of the generating pair so as to ensure that the
exponent sum of $x$ in the relation $r$ is zero. We then have that the
abelianisation of the group $G$ defined by this presentation is equal to
$\z$ if and only if the exponent sum of $y$ in $r$ is $\pm 1$. But if this
is so then $G$ has weight 1 with the normal closure of $x$ being all of $G$.
Moreover, as our group has a deficiency 1 presentation and abelianisation 
$\z$, we obtain $H_2(G;\z)=0$ (for instance by the Hopf formula).

A well known 2-generator 1-relator group that is not residually finite is
the Baumslag-Solitar group 
\[BS(2,3)=\langle x,y|xy^2x^{-1}=y^3\rangle.\]
However another example is one due to Baumslag in \cite{bm} which has
presentation
\[B=\langle x,y|y=y^{-1}x^{-1}y^{-1}xyx^{-1}yx\rangle\]
and this is even less residually finite in that although $B\neq\z$ (in fact it
contains a non-abelian subgroup), every finite quotient of $B$ is abelian
(indeed cyclic as the abelianisation of $B$ is $\z$).
Note that $B$ satisfies the three conditions above in common with classical
knot groups. We can adapt this example to get a group with a badly behaved
commutator subgroup.
\begin{thm}
The group
\[G=\langle a,t|ta^{-2}t^{-1}a^{-1}ta^{-1}t^{-1}atat^{-1}a^{-1}tat^{-1}a
\rangle\]
with abelianisation $G/G'=\z$ has the property that the commutator subgroup
$G'$ is infinite but has no proper subgroups of finite index.
\end{thm}
\begin{proof}
The relation above has exponent sum 0 in $t$ and $-1$ in $a$ so $G'$ is the
kernel of the homomorphism onto $\z$ obtained by setting $a$ equal to the
identity. In order to get a presentation for $G'$ we use the 
Reidemeister-Schreier rewriting process using the obvious Schreier
transversal for $G'$ in $G$ which is $\{t^i:i\in\z\}$. Letting
$a_i=t^iat^{-i}$, we obtain for $G'$ the presentation
\[\langle a_i\mbox{ for }i\in\z|a_{j+1}^{-2}a_j^{-1}a_{j+1}^{-1}a_ja_{j+1}
a_j^{-1}a_{j+1}a_j\mbox{ for }j\in\z\rangle.\]
This means that any two successive generators $a_j,a_{j+1}$ are subject to the
same relation as the one that $x$ and $y$ satisfy in the given presentation
for $B$ (which reveals that the long relation defining
$G$ was obtained by working back from $B$).

Now suppose that $G'$ has a proper finite index subgroup $H$; we can drop down
to a normal subgroup of finite index. We pick any $j$  and restrict the
natural surjective homomorphism $\theta:G'\rightarrow G'/N$ to the subgroup
$\langle a_j,a_{j+1}\rangle$. This subgroup is a quotient of $B$ under
$x\mapsto a_j,y\mapsto a_{j+1}$ but any homomorphism of $B$ into a finite
group must send $y\in B'$ to the identity as the image is abelian.
Thus we have $\theta(a_{j+1})
=\mbox{id}$ and as this is true for any $j\in\z$ we find that $G'=N$, which is
a contraction.
\end{proof}

Note that the group $G$ in Theorem 4.1 is like $B$ in that it too is
an example of a finitely generated group where every finite quotient 
is abelian. This is because if we had a finite index normal subgroup $N$ of
$G$ with $N$ not containing $G'$ then $N\cap G'$ is a proper finite index
subgroup of $G'$. It might be asked whether the converse holds; if
we have a finitely generated
group $G$ such that every finite quotient is abelian then does $G'$ have no 
proper finite index subgroups? In fact we can see this is true in the
case where $G'$ is finitely generated. We argue as follows: if not then
we can take a finite index characteristic subgroup
$C$ of $G'$ which is then normal in $G$, so that $G/C$ is a finitely
generated group which has a finite normal subgroup $G'/C$ and abelian
quotient $G/G'$. This means by \cite{br?} that $G/C$ is virtually abelian
so certainly it is residually finite. But if $G'\neq C$ then  $G/C$ is not
abelian so on taking an element $\gamma$ in the commutator subgroup of
$G/C$ and a normal finite index subgroup of $G/C$
missing $\gamma$, we get a finite
non-abelian quotient of $G/C$ and hence of $G$. However 
we can use K.\,S.\,Brown's algorithm in \cite{bks} to show that the commutator
subgroup of $B$ and of $G$ in Theorem 4.1 are both not finitely generated.

These group theoretic constructions are all very well but here we can realise
them geometrically. As already mentioned, the group $G$ in Theorem 4.1 cannot
be a classical knot group because it is not residually finite. However we
can move up to higher dimensions and consider an $n$-knot, that is a locally
flat embedding of $S^n$ into $S^{n+2}$, and the fundamental group of its
complement (see \cite {hill} for details).
Necessary conditions that the fundamental group of this
complement must possess are the three given above for classical knot groups,
and it was shown by Kervaire in \cite{k65bsmf} that these conditions are
sufficient when $n\geq 3$, whereupon we have 
a smoothly embedded $S^n$ in $S^{n+2}$.
This is not true for $n=2$ (\cite{lev77}) but if we 
strengthen the condition $H_2(G;\z)=0$ to $G$ has
deficiency 1 then we do now have sufficient conditions for $n=2$ (though now
these are no longer necessary). This was also proved by Kervaire 
in \cite{k65bsmf}, at least on
assumption of the topological 4 dimensional Poincare conjecture that was
later established by Freedman \cite{fr}. Furthermore if we take the
presentation of the trivial group given by setting $t$ equal to the identity,
where $t$ is the element whose normal closure is all of $G$, then by 
\cite{lev77} we can
assume that we have a smooth 2-knot in the standard 4-sphere if this
presentation can be made trivial by Andrews-Curtis moves. This allows us to
express Theorem 4.1 in purely geometric terms.

\begin{co} 
There exists a smooth knotted 2-sphere in $S^4$ such that the only finite
covers of its complement are cyclic and such that its infinite cyclic cover
is not simply connected but has no proper finite covers at all.
\end{co}
\begin{proof}
Given any 2-generator 1-relator group $\langle x,y|r\rangle$ where the
exponent sum of $x$ in $r$ is zero and the exponent sum of $y$ in $r$ is
$\pm 1$, we have that $\langle x,y,|r,x\rangle$ is trivialisable using
Andrews-Curtis moves because we can move appearances of $x$ from the
beginning and end of $r$, and then we can conjugate by appropriate powers
of $y$ to put $x^{\pm 1}$ on the beginning or end of the relator again. 
This will terminate in a word of the form $y^k$ but the moves do not change
the exponent sum of $y$ at each stage. Thus
the group $G$ in Theorem 4.1 is the fundamental group of the complement of
a smooth 2-sphere in $S^4$, with the infinite cyclic cover having
fundamental group $G'$. If there were a finite index subgroup of $G$ which
did not contain $G'$ then its intersection with $G'$ would be a proper finite
index subgroup of $G'$.
\end{proof}

We also get as an immediate corollary from Kervaire's characterisation of high
dimensional knot groups and Theorem 4.1 that there exists a 
(smoothly) knotted $n$-sphere
in $S^{n+2}$ for all $n\geq 2$ such that the  only finite
covers of its complement are cyclic and such that its infinite cyclic cover
is not simply connected but has no proper finite covers. 

We finish by mentioning an intriguing conjecture which is Conjecture 4.4
in \cite{swtk} that if $k$
is a non-fibred classical knot then there exists a finite group $F$ such
that the dynamical shift space of representations from the commutator
subgroup of $\pi_1(k)$ into $F$ has strictly positive topological entropy.
This would then distinguish fibred knots from non-fibred knots
because if the commutator subgroup
is finitely generated then the space of representations into any finite
group is finite and hence the entropy is zero.
However we see from the above examples that the
conjecture fails if we enlarge it to claim that it determines whether or not
a finitely presented group with a homomorphism $\chi$ onto $\z$ has its
kernel finitely generated, or even that the complement of an $n$-knot in
$S^{n+2}$ is fibred, because in Theorem 4.1 and Corollary 4.2 the
commutator subgroup of $G$ is infinitely generated but the space of
representations of $G'$ into any finite group is trivial. We also note that
Problem 2.6 in this paper asks for which homomorphisms $\chi$ of
finitely presented groups onto $\z$ do we find that the dynamical shift
space of representations from $\mbox{ker }\chi$ to any finite group is
finite. There 
Proposition 2.7 states that this happens if $\mbox{ker }\chi$ does
not contain a free group of rank 2, but we see that the converse is not true
because in Theorem 4.1 $G$ and hence $G'$ contain non-abelian free subgroups.

\end{document}